\documentclass[microtype]{gtpart}
\agtart
\usepackage{amscd, ragged2e}

\makeatletter
\def\thebibliography#1 {\@thebibliography@{99999}\small\parskip0pt 
plus2pt\relax\addcontentsline{toc}{section}{Bibliography}}
\makeatother

\title{Quantum invariants of random 3-manifolds}

\author{Nathan M.~Dunfield}
\address{Dept.~of Mathematics, MC-382, University of Illinois, Urbana, IL 61801, USA}
\email{nathan@dunfield.info}
\urladdr{http://dunfield.info}

\author{Helen Wong} \address{Dept.~of Mathematics, 1 North College St., Northfield, MN 55057, USA}
\email{hwong@carleton.edu}
\urladdr{http://people.carleton.edu/~hwong}

\begin{abstract}
  We consider the $\SO 3$ Witten-Reshetikhin-Turaev quantum invariants
  of random 3-manifolds.  When the level $r$ is prime, we show that
  the asymptotic distribution of the absolute value of these
  invariants is given by a Rayleigh distribution which is
  independent of the choice of level.  Hence the probability that the
  quantum invariant certifies the Heegaard genus of a random
  3-manifold of a fixed Heegaard genus $g$ is positive but very small,
  less than $10^{-7}$ except when $g \leq 3$.  We also examine random
  surface bundles over the circle and find the same distribution for
  quantum invariants there.
 \end{abstract}


\arxivreference{arXiv:1009.1653}
\arxivpassword{63dby}

\volumenumber{}
\issuenumber{}
\publicationyear{}
\papernumber{}
\startpage{}
\endpage{}
\doi{}
\MR{}
\Zbl{}
\received{}
\revised{}
\accepted{}
\published{}
\publishedonline{}
\proposed{}
\seconded{}
\corresponding{}
\editor{}
\version{}

\def\RCS$#1: #2 ${\expandafter\def\csname RCS#1\endcsname{#2}}
\RCS$Revision: 62 $
\RCS$Date: 2011-07-03 08:54:56 -0500 (Sun, 03 Jul 2011) $
\renewcommand{\today}{\number\year /\number\month /\number\day}
\newcommand{\versioninfo}{Version \RCSRevision; Last commit \RCSDate; Compile \today }

\newcommand{\F}{{\mathbb F}}
\newcommand{\maps}{\colon\thinspace}
\DeclareMathOperator{\tr}{tr}
\newcommand{\SU}[1]{ {\mathrm{SU}({#1})}}
\newcommand{\SO}[1]{ {\mathrm{SO}({#1})}}
\newcommand{\U}[1]{ {\mathrm{U}({#1})}}
\newcommand{\PU}[1]{ {\mathrm{PU}({#1})}}
\newcommand{\rightquom}[4]{{\raisebox{#3}{$#1$}}#4/ {\raisebox{-#3}{$#2$}}}
\newcommand{\MCG}{\Gamma}
\newcommand{\MCGext}{\tilde \Gamma_1}
\newcommand{\Dns}{\mathcal{D}_{ns}}

\newcommand{\talpha}{t_{\alpha}}

\newcommand{\labelset}{\mathcal{L}}

\newcommand{\firstroot}{ e ^{2 \pi i / 4r}}
\newcommand{\Free}{\mathrm{Free}}

\theoremstyle{plain}

\swapnumbers
\newtheorem{theorem}[subsection]{Theorem}
\newtheorem{conjecture}[subsection]{Conjecture}
\newtheorem{lemma}[subsection]{Lemma}

\newtheorem{proposition}[subsection]{Proposition}

\theoremstyle{definition}

\theoremstyle{remark}
\newtheorem{remark}[subsection]{Remark}

\theoremstyle{plain}
\newtheorem{mainonebackwards}{Theorem}

\makeatletter
  \let\c@figure=\c@subsection
  
  \let\c@table=\c@subsection
  
  \let\c@equation=\c@subsection
  \let\theequation=\thesubsection
\makeatother

\newcommand{\hyp}{\nobreakdash-\hspace{0pt}}
\newcommand{\3}[1]{3\hyp}

\begin{document}
\maketitle

\section{Introduction}
\label{sec_intro}

An important class of \3-manifold invariants arises from the
topological quantum field theories (TQFTs) introduced by Witten
\cite{Witten1989} and then formalized by Atiyah \cite{Atiyah1990a,
  Atiyah1990b}.  At the coarsest level, a TQFT associates to each
closed \3-manifold $M$ a complex number $Z(M)$.  Here, we study the
distribution of $Z(M)$ for a random \3-manifold in the sense of
Dunfield and W.~Thurston \cite{DunfieldThurston:random}.  Such
manifolds are generated by gluing together two handlebodies of large
genus by an extremely complicated element in the mapping class group
(see Section~\ref{sec-random} for details).  For the $\SO 3$ quantum
invariants at a prime level, we show that the distribution of $|Z(M)|$
is in fact independent of the particular invariant:
\begin{theorem} \label{Thm:DistributionAbs} Let $Z$ be the $\SO 3$
  TQFT associated to a prime level $r \geq 5$.  Then for a random
  \3-manifold $M$, the invariant $|Z(M)|$ is distributed by the
  Rayleigh distribution with mean $\sqrt{\pi} / 2$.  In particular,
  the probability that $|Z(M)|$ is at least some $x \in [0, \infty)$
  is:
  \[
  P\big\{ |Z(M)|  \geq  x \big\} = e^{-x^2} 
  \]
\end{theorem}
\noindent
We conjecture that $Z(M)$ itself is distributed by the standard
complex Gaussian, though we were not able to prove this.  (See
Section~\ref{Section:Distribution} for the definitions of
these probability distributions.)  It is known
that for each prime level the values of $Z(M)$ form a dense set
in $\C$ \cite{Wong:density}.

While Theorem~\ref{Thm:DistributionAbs} has intrinsic interest, our
motivation was a classical topological question: when does the rank of
$\pi_1(M)$ determine the \emph{Heegaard genus} of $M$, that is, the
minimal genus of a Heegaard splitting of $M$? This was first posed by
Waldhausen \cite{Haken1970, Waldhausen1978}, and in 1984 Boileau and
Zieschang gave the first breakthrough examples where the rank and
genus differ \cite{BoileauZieschang1984}.  However, when $M$ is
hyperbolic, there are still no known examples where the rank and genus
differ (see \cite{Souto2007} for a survey).  Garoufalidis
\cite{Garoufalidis1998} and Turaev \cite{Turaev1994} showed that
$|Z(M)|$ gives a lower bound on Heegaard genus (see
Theorem~\ref{Thm:HGbound} below), and Wong used this to give an
alternate proof that the examples of Boileau and Zieschang have rank
smaller than genus \cite{Wong:BZ}.  This is interesting in part
because the technique of \cite{BoileauZieschang1984} cannot be applied
to any hyperbolic 3-manifold, whereas quantum invariants have no such
apparent restrictions.  This led us to try to find hyperbolic
3-manifolds where the rank and genus differ and the genus could be
certified by $|Z(M)|$.

Unfortunately, in that search we found that quantum invariants
typically gave very poor bounds on the Heegaard genus.  This can be
explained by our results here.  In particular,
Theorem~\ref{Lem:ZMdistributedasXd} gives the precise distribution of
$|Z(M)|$ for a random \3-manifold of Heegaard genus $g$, and we use
this to show:
\begin{theorem}\label{thm-genus-bound} 
  Fix a prime level $r \geq 5$.  Let $M$ come from a random Heegaard
  splitting of genus $g$.  Then the $\SO{3}$ quantum invariant of
  level $r$ gives a sharp lower bound on the genus of $M$ with
  probability $\left(1 - \mu^2 \right)^{d -1}$, where $\mu =
  \frac{2}{\sqrt{r}} \sin \frac{\pi}{r}$, and $d$ is the dimension of
  the Hilbert space of this TQFT for a surface of genus $g$.
\end{theorem}

Here $d$ depends on both $g$ and $r$ with $d \sim \mu^{2 - 2 g}$ (see
Lemma~\ref{Lem:MuVersusD}); since $\mu < 1$, the dimension $d$
increases rapidly in both $r$ and $g$.  Table~\ref{tab-genus-bound}
gives the probabilities when $r$ and $g$ are small.
\begin{table}[thb]
  \begin{center}
    \renewcommand{\arraystretch}{1.3}
    \addtolength{\tabcolsep}{-2pt}
    \begin{tabular}{l|ccccc}
      & $r = 5$ & $r = 7$ & $r = 11$ & $r = 13$ & $|H_1(M;\F_2)|$\\ \hline
      $g = 2$ & $0.274$ & $0.228$ & $0.206$ &$0.202$ & $0.067$ \\
      $g = 3$ & $0.011$ & $1.6 \times 10^{-5}$& $1.2 \times 10^{-17}$& $4.0\times10^{-28}$ & $0.007$\\
      $g = 4$ & $1.3\times10^{-7}$ & $7.5 \times 10^{-42} $ & $1.8 \times 10^{-542}$ & $1.2 \times 10^{-1442} $& $4.3 \times 10^{-4}$\\
      $g = 5$ & $3.6\times10^{-25}$& $1.2\times10^{-373}$& $4.2 \times 10^{-18437}$ & $6.7 \times 10^{-80495}$& $1.3 \times 10^{-5}$\\
    \end{tabular} 
  \end{center}
  
  \vspace{0.3cm}
  \caption{Probabilities
    that the $\SO{3}$ quantum invariants of various levels $r$ give
    sharp bounds on genus.  For comparison, the far right
    column gives the corresponding probability that the rank of $H_1(M;
    \F_2)$ certifies the genus, namely $\prod_{k = 1}^g (2^k + 1)^{-1}$ as determined by 
    \cite[\S8.6]{DunfieldThurston:random}.  }\label{tab-genus-bound}
\end{table} 
As you can see, the probability that $Z(M)$ gives a sharp genus bound
decreases as either $g$ or $r$ increases, and the bound is very rarely
sharp except when $g = 2$ or (marginally) when $(g,r) = (3,5)$.
Unfortunately, the case $g=2$ isn't interesting for the question of
rank versus genus since by Perelman's Geometrization Theorem, any
closed \3-manifold where $\pi_1(M)$ is cyclic must have genus one.
Moreover, it seems more likely that hyperbolic examples with rank differing
from genus exist when $g > 2$.  However, Table~\ref{tab-genus-bound}
makes it clear that it would be very hard to find such examples where
the genus is certified by $Z(M)$ via a brute force search.  Of course,
it is possible to systematically generate examples where $Z(M)$ does
certify the genus by always doing certain powers of Dehn twists, but
we had no luck finding examples among such manifolds.

One natural class of hyperbolic \3-manifolds is surface bundles over
the circle, and one place we looked for possible examples was the
Hall-Schleimer census of small genus 2 bundles \cite{HallSchliemer}.
For the natural notion of random surface bundle, we show here that the
asymptotic distribution of $|Z(M)|$ is the same as that for all random
\3-manifolds, even though general random \3-manifolds are bundles with
probability 0.  This can be taken as evidence for the naturality and
robustness of this notion of random \3-manifold.  We also show that
for a generic bundle, no $\SO{3}$ quantum invariant certifies its
Heegaard genus.  For details, see Section~\ref{sec-bundle}.

\subsection{Outline}

In Section~\ref{sec-random}, we review the notion of a random
\3-manifold and make precise how we discuss (asymptotically) their
various properties.  Then in Section~\ref{Section:CentralExt}, we
introduce our preferred formalism for $Z(M)$, in particular focusing
on certain (projective) unitary representations of mapping class
groups.  Section~\ref{Section:Distribution} is the core of the paper
and contains the proofs of the main results.  The two main tools we
use are the work of Freedman, Larsen, and Wang
\cite{LarsenWang2005, FreedmanLarsenWang2002} on the density of these
unitary representations, and the Ito-Kawada theorem about
equidistribution of random walks in compact groups.  Finally,
Section~\ref{sec-bundle} discusses what happens for surface bundles.

\subsection{Acknowledgments}  

Dunfield was partially supported by US NSF grant DMS-0707136 and a
Sloan Fellowship.  Some of this work was done while Wong was at
Bowdoin College.  The authors thank Zhenghan Wang for his suggestion
to consider this problem and other helpful discussions, and also
  thank the anonymous referee for their comments. 

\section{Random Heegaard splittings}
\label{sec-random}

We now review the notion of \emph{random \3-manifold} studied in
\cite{DunfieldThurston:random}.  Fix a positive integer $g$.  Let $H$
be a genus-$g$ handlebody, and denote $\partial H$ by $\Sigma$.  Let $\Gamma$ be the
mapping class group of $\Sigma$.  For any mapping class $f \in \Gamma$, we may
associate to it a closed \3-manifold $M_f$ obtained by gluing
together two copies of $H$ via the homeomorphism $f$.

To define a random Heegaard splitting of genus $g$, fix generators $T$
for $\Gamma$.  A random element $w$ of $\Gamma$ of complexity $\ell$ is defined to
be the result of a random walk in the generators $T$ of length $\ell$,
that is, a random word of length $\ell$ in $T$.  For such a $w$, we call
$M_{w}$ the \emph{manifold of a random Heegaard splitting} of
genus $g$ and complexity $\ell$. 

We are interested in the properties of such random $M_w$ as the
complexity $\ell \to \infty$.  A priori, this might depend on the choice of
generators $T$ for $\Gamma$, though it turns out that many properties
behave nicely and are independent of such choices.  Consider a
function
\[
F \maps (\mbox{closed \3-manifolds}) \to S,
\]
where $S$ is some set; a typical example is $F(M) = \dim H_1(M;
\F_2)$.  Now let $\mu_\ell$ be the finitely supported atomic measure on $S$
defined by taking $\mu_\ell(\{s\})$ to be the probability that $F(M) = s$,
where $M$ is a random Heegaard splitting of genus $g$ and complexity
$\ell$.  If there is a measure $\nu_g$ on $S$ so that, for every choice of
generating set of $\Gamma$, the measures $\mu_\ell $ converge weakly to $\nu_g$
as $\ell \to \infty$, then we say that \emph{$F$ is distributed by $\nu_g$ for a
  random \3-manifold of genus $g$}.  When $\nu_g$ is atomic, we will say
things like ``the probability that $F(M) = s$ is $\nu_g(s)$'', though
properly speaking this is a limiting statement as $\ell \to \infty$.

If the $\nu_g$ themselves have a weak limit $\nu$ as $g \to \infty$, then we say that
\emph{$F$ is distributed by $\nu$ for a random 3-manifold}.  Of
course, not all $F$ have such limiting distributions $\nu$, but
for instance $F(M) = \dim H_1(M;\F_p)$ does
\cite[\S8.7]{DunfieldThurston:random}.  Another example is that, for
fixed $n$, the number of epimorphisms from $\pi_1(M)$ to the
alternating group $A_n$ is Poisson distributed with mean 1
\cite[Thm~7.1]{DunfieldThurston:random}.  A geometric example is
Maher's result that $M$ is hyperbolic with probability 1 \cite{Maher2008}.

\section{WRT invariants via unitary representations} \label{Section:CentralExt}

In this section, we review the definition of the
Witten-Reshetikhin-Turaev TQFT for $\SO{3}$ in terms of a certain
unitary representation of an extension of the mapping class group.
Though the WRT TQFT was originally defined using the representation
theory of $U_q\mathit{sl}(2)$ in \cite{Witten1989} and \cite{ReshTuraev1991},
an equivalent definition using skein theory was developed in
\cite{BHMV1995} and \cite{Lickorish1991}, among others.  We follow the
latter approach as summarized in \cite{Lickorish2002}, and will
also refer to \cite{MasbaumRoberts1995} when necessary.

Fix an odd integer $r \geq 5$ as the level, and let $A = i\firstroot$.
For a closed surface $\Sigma$ of genus $g$, the associated
finite-dimensional Hilbert space $(V, \langle \cdot, \cdot \rangle)$
has the following form.  Here $V$ depends on both $r$ and $g$, and so
we will also denote it by $V_{r,g}$ when this needs to be
emphasized. Let $\labelset = \{ 0, 2, 4, \ldots r-3 \}$ be the set of
even labels.  The basis elements of $V$ correspond to the admissible
$\labelset$-labelings of trivalent spines of a handlebody $H$ with
$\partial H = \Sigma$.  The dimension of $V$ is given by the Verlinde
formula \cite{BHMV1995} as
\begin{equation}\label{eq-verlinde}
  d =   \frac{1}{2^g} \left(\frac{r}{2}\right)^{g-1} \sum_{j=1}^{r-1} \left(\sin \frac{\pi j}{r} \right)^{2-2g}  
\end{equation}
Moreover, for our choice of $A$, there is an inner product $\langle
\cdot, \cdot \rangle$ induced from the splitting $\#_g S^1 \times S^2
= H \cup (-H)$, and these basis elements are orthogonal with respect
to $\langle \cdot, \cdot \rangle$.  See \cite{Lickorish2002} for
details.

Underlying this inner product, and indeed the entire WRT TQFT theory,
is a skein element $\omega$ which has the key property that if two
links related by Kirby handleslides are colored by $\omega$, then the two
resulting skeins have the same Kauffman bracket evaluation.  More
precisely, we have chosen $\omega$ so that the Kauffman bracket of a
$0$-framed unlink in $S^3$ colored by $\omega$ is the number $\mu^{-1}$,
where $\mu^2 = -\big(A^2 - A^{-2}\big)^2/r$ and so $\mu = \frac{2}{\sqrt{r}} \sin
\frac{\pi}{r}$.  We denote the Kauffman bracket of a $(-1)$-framed
unlink colored by $\omega$ by $\kappa^3$.  Formulas for $\kappa^3$ appear in
e.g.~\cite{BHMV1995}, but here all that is needed is that $\kappa^3$ is a
$4r^{\mathrm{th}}$ root of unity.

We now describe the (projective) action of the mapping class group
$\Gamma$ of $\Sigma$ on $V$.  Let $\Dns$ denote the set of Dehn twists
along non-separating curves on $\Sigma$.  These generate $\Gamma$,
and to start, we define an action of the free group $\Free(\Dns)$ on $V$.
To a word $w$ in $\Free(\Dns)$, we associate a link $L(w)$ whose
components have framing $\pm 1$ and that lives in $\Sigma \times I$,
as in the proof of the Dehn-Lickorish theorem.  Color the components
of $L(w)$ by the skein element $\omega$ to obtain $s(w)$.  Consider
the transformation on $V$ which adds to a collar neighborhood of $H$
the skein $s(w)$.  Because $\omega$ behaves well under handleslides,
this gives a unitary transformation on $V$.  In other words, we have a
homomorphism $\rho \maps \Free(\Dns) \to \U{V}$.  This construction is
the ``skein-theoretic version of the original geometric action on Dehn
twists" from \cite{MasbaumRoberts1995}.  While the representation of
$\Free(\Dns)$ on $V$ does not descend to $\Gamma$ itself, it does give
a \emph{projective} representation of $\Gamma$, i.e.~there is a
homomorphism $\overline \rho: \Gamma \to \PU{V}$ with the property
that
\begin{equation} \label{CD}
 \begin{CD}
 \Free(\Dns) @>  \rho >>  \U{V}\\
@VVV @VV V\\
\MCG @> \overline \rho >>  \PU{V}
\end{CD}
\end{equation}
commutes.  

The WRT 3-manifold invariant may be described easily using this
language of unitary representations.  For any word $w$ in
$\Free(\Dns)$, let $M_w$ be the closed \3-manifold obtained from the
corresponding element of $\Gamma$.  The \emph{vacuum vector} $v_\emptyset \in
V$ is the element corresponding to the empty link.  Up to a phase
factor involving $\kappa^3$, the action of $\rho(w)$ on $v_\emptyset$ determines
$Z(M_w)$.  The phase factor is in some sense a correction factor for
the framing, and depends on the exponent sum $e(w)$ of the word $w$
and on the signature $\sigma_b(w)$ of the linking matrix of $L(w)$.
Precisely, the WRT invariant is given by
\begin{equation}
Z(M_w) =   \kappa ^{-3( \sigma_b(w) + e(w))}  \cdot \mu \cdot \langle \rho(w) v_\emptyset, v_\emptyset \rangle
\end{equation}
For example, since the vacuum vector has $ | v_\emptyset|^2 = \mu^{-g}$, we have
$Z(S^3) = \mu$ and $Z(S^1 \times S^2) = 1$; in general $Z(\#_g S^1 \times S^2) =
\mu^{1-g}$.  Note that $Z(M)$ here differs from $I_p(M)$ appearing in
\cite{MasbaumRoberts1995} by a factor of $\mu$.

It will be convenient to rewrite $Z(M_w)$ using matrix notation.  Take
the normalized vacuum vector $\frac{v_\emptyset}{ | v_\emptyset|}$ as the first basis
element in an orthonormal basis for $V$, and denote the $(1,1)$-entry
of $V(w)$ by $\rho(w)_{(1,1)}$.  Then,
\begin{equation}\label{eq-Z-mat}
Z(M_w) = \kappa^{-3 (\sigma_b(w) + e(w))} \cdot \mu^{1-g} \cdot \rho( w)_{(1,1)},
\end{equation}
and since $\kappa$ is a root of unity and $\mu$ positive, we have
\begin{equation}\label{eq-abs-Z-mat}
|Z(M_w)| = \mu^{1-g} | \rho( w)_{(1,1)}|.
\end{equation}
In particular, $|Z(M_w)| \leq \mu^{1 - g}$, leading to the key
observation of Garoufalidis \cite{Garoufalidis1998} and Turaev
\cite{Turaev1994}:
\begin{theorem}\label{Thm:HGbound}
  Suppose $M$ is a closed orientable 3-manifold with Heegaard
  genus $g(M)$.  Then $ |Z(M)| \leq \mu^{1-g(M)}$.
\end{theorem}

\subsection{The image of $\rho$}

We now turn to the properties of $\rho$ that will be needed later.    The first proposition asserts that the image of $\rho$ is contained in $\SU{V}$, or in some exceptional cases, a slightly larger group.    

\begin{proposition} \label{Prop:SU} Let $r \geq 5$ be odd. The image of
  the homomorphism 
  \[
  \det \circ \rho \maps \Free(\Dns) \to \C^\times 
  \] 
  is $1$ for all $g \geq 2$, except when both $g=2$ and $r \equiv 0
  \pmod{5}$.  In the latter case, or when $g = 1$, the image is
  contained in the subgroup generated by $A^4$.  
\end{proposition}

\begin{proof}
  When $g \geq 3$, we need to show that $\det \circ \rho \maps \Free(\Dns) \to
  \C^\times$ is the trivial homomorphism.  Following
  \cite{MasbaumRoberts1995}, there is a central extension $\MCGext$ of
  $\MCG$ to which $\rho$ descends, and in this proof we consider
  $\MCGext$ as the domain of $\rho$.  As discussed in the proof of
  Theorem 3.10 of \cite{MasbaumRoberts1995}, the abelianization
  $\MCGext^{\mathit{ab}} = H_1(\MCGext; \Z)$ is $0$ when $g \geq 3$.
  Thus the homomorphism $\det \circ \rho \maps \MCGext \to \C^\times$ must be
  trivial as claimed, since it factors through the trivial group.
 
  So now we focus on $g = 1, 2$, where $\MCGext^{\mathit{ab}} =
  H_1(\MCGext; \Z) = \Z$.  By Theorem~3.8 and Remark~3.9(ii) of
  \cite{MasbaumRoberts1995}, all the Dehn twists in $\Dns$ become
  conjugate in $\MCGext$.  Thus the abelianization $H_1(\MCGext; \Z)$
  is generated by any single non-trivial element $\talpha$ of $\Dns$.
  We proceed by calculating the determinant of $\rho(\talpha)$
  explicitly, and find $\rho(\talpha) = A^e$ where 
  \[ 
  e =  \begin{cases} {\frac{1}{6} (r-3)(r-1)(r+1)} & \mbox{when $g= 1$,}\\
    {\frac{r}{5} \big(\frac{r-3}{2}\big)\big(\frac{r-1}{2}\big)\big( \frac{r+1}{2}\big)\big(\frac{r+3}{2}\big)}
    & \mbox{when both $g=2$ and $r \equiv 0 \bmod{5}$, }\\
    0 & \mbox{in all other cases.}
  \end{cases}
  \] 
  In the first two cases, it is easy to see that $4 | e$, and so the
  image of $\det \circ \rho$ is contained in $\left\langle A^4\right\rangle \leq \C^\times$.

  First, consider the case $g=2$.  We view the genus 2 handlebody $H$ as a neighborhood of the theta
  graph, with edges labeled by $ a,b,c \in \{0, 2, 4, \ldots , r-3\}$.  The
  basis elements of $V$ correspond to those labelings that are
  admissible, namely those satisfying
  \[
  a+b+c \leq 2(r-2), \quad \quad b+c \geq a, \quad \quad a+ c \geq b, \quad
  \quad a+ b \geq c.
  \] 
  Consider a simple closed curve $\alpha \subset \Sigma$ which bounds a disc in $H$
  which is dual to the edge labeled by $a$.  Then each basis element
  $(a,b,c)$ of $V$ is an eigenvector for $\rho(\talpha)$, with
  eigenvalue $A^{a^2 + 2a}$.  Now $\det \rho(\talpha)$ is the product
  of these eigenvalues, and hence is $A^e$ for $e = \sum_{(a,b,c) \,
    \mbox{\scriptsize admissible}} a^2 + 2a$.

  For a fixed label $a$, it can be shown that the number of admissible
  triples containing $a$ is exactly $(a+1)( \frac{r-1- a}{2} )$.  Thus
  we can apply standard formulas to show that
  \begin{eqnarray*}
  e &=&  \sum_{\stackrel{(a,b,c)}{\mbox{\tiny admissible}}} a^2 + 2a
   = \sum_{\mbox{\tiny even } a=0}^{r-3} (a^2 + 2a) (a+1)\left( \frac{r-1- a}{2} \right)\\
   & = & \sum_{k=0}^{ \frac{r-3}{2}}  -8k^4 + (4r-16) k^3 + (6r-10)k^2 + (2r-21) k\\
   & = & \frac{r}{5}  \left(\frac{r-3}{2}\right)\left( \frac{r-1}{2}\right)\left( \frac{r+1}{2}\right)\left( \frac{r+3}{2}\right).
 \end{eqnarray*}
 To conclude, we claim $4r$ divides $e$ and hence $\det \rho(\talpha)
 = A^e = 1$ since $A$ is a $4r^{\mathrm{th}}$ root of unity.  As $5
 \nmid r$, then 5 must divide one of the consecutive
 integers $\frac{r-3}{2}$, $\frac{r-1}{2}$, $\frac{r+1}{2}$,
 $\frac{r+3}{2}$ and hence $r | e$.  It is also clear that $4$ divides
 $a^2 + 2a$ since $a$ is even, and hence $4 | e$.  Thus $4r | e$ as
 needed.
 
 When $g=1$, we view the solid torus $H$ as a neighborhood of a
 circle, with the Dehn twist curve $\alpha$ being a meridian linking the
 circle once.  Thus the eigenvalues of $ \rho(\talpha)$ are $A^{a^2 +
   2a}$ all with multiplicity one.  Then $e = \sum_{\mbox{\tiny even }a =
   0}^{r-2} a^2 + 2a = \frac{1}{6} (r-3)(r-1)(r+1)$, as claimed.
\end{proof}

The next result is a simple corollary of the results of Freedman, Larsen,
and Wang \cite{FreedmanLarsenWang2002, LarsenWang2005} that the image
of $\overline \rho$ is dense in $\PU{V}$ when $r$ is an odd prime.
\begin{proposition} \label{Prop:DenseImage} Let $G \leq \U{V}$ be the
  closure of the image of $\rho$.  If $g \geq 2$ and $r \geq 5$
    is prime,
  the group $G$ contains $\SU{V}$.  Moreover $G = \SU{V}$ except if
  both $g=2$ and $r = 5$.
\end{proposition}
\begin{proof}
  First suppose that $(g,r) \neq (2,5)$, so that $G \subset \SU{V}$ by
  Proposition~\ref{Prop:SU}.  Since $\rho$ is a homomorphism, $G$ is a
  compact Lie subgroup of $\SU{V}$ and hence $\dim G \leq \dim \SU{V}$.
  On the other hand, the results of Freedman, Larsen, and Wang
  (Theorem 3 of \cite{LarsenWang2005}) show that the map $\overline
  \rho$ has dense image.  Thus the commutativity of (\ref{CD}) implies
  that $G$ maps onto $\PU{V}$, and so $\dim G \geq \dim \PU{V}$.  As
  $\dim \SU{V} = \dim \PU{V}$ and $\SU{V}$ is connected, we must have
  $G = \SU{V}$ as claimed.

  If $(g,r) = (2,5)$, then we saw in the proof of
  Proposition~\ref{Prop:SU} that the image of $\det \circ \rho$ is
  contained in the finite subgroup of $\C^*$ generated by $A^4$.  Let $G_0$
  be the connected component of the identity, which must be contained
  in $\ker ( \det \circ \rho )$ since the image of $\det \circ \rho$ is finite.
  Now as before, we must have $\dim G \geq \dim \PU{V}$.  Since $G_0 \leq
  \SU{V}$ and $\dim G_0 = \dim G$ this forces $G_0 = \SU{V}$ as
  needed.
\end{proof}

\begin{remark}\label{rem-denseimage}
  Note that Proposition~\ref{Prop:DenseImage} remains true if we
  restrict the domain of $\rho$ to any subgroup of $\Free(\Dns)$ whose
  image generates $\Gamma$; this is because the proof only uses that the
  image of $\overline \rho$ is dense $\PU{V}$.
\end{remark}

We will also need to understand the dimension of $V = V_{r,g}$ as a function of
$r$ and $g$ (recall here that $\mu$ depends on $r$).  
\begin{lemma} \label{Lem:MuVersusD} Fix an odd integer level $r \geq 5$.
  Let $d_g$ be the dimension of the WRT representation space $V_{r,g}$.
  Then \[ \lim_{g \to \infty} \frac{d_g}{\mu^{2-2g}} = 1. \] As $\mu <
  1$, this implies $d_g \to \infty$ as $g \to \infty$.
\end{lemma}
\begin{proof}
  Using (\ref{eq-verlinde}) and the formula for $\mu$ we have:
  \[
  \frac{d_g}{\mu^{2 - 2g}} = \frac{1}{2} \sum_{j=1}^{r-1}
  \left(\frac{\sin(\pi j/r)}{\sin(\pi/r)} \right)^{2 - 2g} =
  \frac{1}{2} \sum_{j=1}^{r-1} a_j^{2g-2} \quad \mbox{where} \quad a_j = \frac{\sin(\pi/r)}{\sin(\pi j/r)}
  \]
  Now when $j= 1$ or $r-1$, then $ \sin ( \pi j / r) = \sin(\pi/r)$
  and so $a_j = 1$.  However, for any other $j$, we have $ \sin(\pi j
  / r) > \sin(\pi/r)$ and so $0 \leq a_j < 1$.  Thus
  \begin{equation}\label{eq-dg-vs-mu}
  \frac{d_g}{\mu^{2 - 2g}} = \frac{1}{2} \left(2 + \sum_{j=2}^{r-2}
    a_j^{2g - 2}\right)
  \end{equation}
  where all $a_j$ left are $<1$ and hence $a_j^{2g - 2} \to 0$ as $g
  \to \infty$.  As $r$ is fixed, we see that $\lim_{g \to \infty}
  {d_g/\mu^{2 - 2g}} = 1$ as claimed.  Indeed, it converges 
  exponentially fast in $g$. 
\end{proof}

\section{Distribution of $|Z(M)|$} \label{Section:Distribution}

In this section, we calculate the distribution of $|Z(M)|$ for random
\3-manifolds.  To state our results, we first define some
distributions.  The \emph{standard complex Gaussian} is the
probability distribution on $\C$ where the real and imaginary
components of $z$ are independently distributed Gaussians of mean 0
and variance $1/2$; equivalently the distribution function is
$\frac{1}{\pi}e^{-|z|^2}$ (see e.g.~\cite[page 5]{Novak2007}) and
hence the expectation of $|z|^2$ is $1$.  The probability distribution
describing of the absolute value of a standard complex Gaussian is the
Rayleigh distribution with mean $\frac{\sqrt{\pi}}{2}$; explicitly the
density function is $2 x e^{-x^2}$ and $P\{ R \geq x \} = e^{-x^2}$
(see e.g.~\cite{Papoulis2002}).  We refer to this distribution here as
the \emph{standard Rayleigh distribution}; however, some authors use
this term for the Rayleigh distribution with mean $\sqrt{\pi/2}$ which
arises when the components of the original complex distribution have
variance $1$. Recall our main result is:
\begin{mainonebackwards}
  Let $Z$ be the $\SO 3$ TQFT associated to a prime
  level $r \geq 5$.  Then for a random \3-manifold $M$, the invariant
  $|Z(M)|$ is distributed by the standard Rayleigh.
\end{mainonebackwards}

Along the way we give a complete picture for each fixed Heegaard
genus, using the following distribution.  Let $X_d$ be the random
variable on $\C$ which gives the distribution of the (1,1)-entry of a
unitary matrix of size $d$ with respect to Haar measure; it is
well-known, see e.g.~\cite[pg.~140]{HiaiPetz} that the density in
polar coordinates is $\frac{d-1}{\pi}(1 - r^2)^{d - 2} \, \mathit{dA}$
for $0 \leq r \leq 1$.  Focusing on the absolute value $|X_d|$, we easily calculate
\[
P\left\{ |X_d| \geq x \right\} = \big(1 - x^2\big)^{d - 1}  \quad \mbox{ for  } \,  0 \leq x \leq 1. \]
We first show
\begin{theorem} \label{Lem:ZMdistributedasXd} Fix a genus $g$ and
  prime level $r \geq 5$.  Let $d_g$ be the dimension of the Hilbert
  space $V_{r,g}$ associated to the surface of genus $g$.   Then the invariant $|Z(M)|$ is distributed by
  $\mu^{1-g} |X_{d_g}|$ for a random Heegaard splitting of genus $g$.  In particular
  \[
  P\left\{ |Z(M)| \geq x \right\} = \big(1 - \mu^{2 g - 2}x^2\big)^{d_g - 1}  \quad \mbox{ for }\,  0 \leq x \leq \mu^{1 - g}.
  \]
\end{theorem}

\begin{proof}
  Consider a set $T$ consisting of words $t_1, \ldots t_n \in \Free(\Dns)$
  which together generate $\Gamma$.  Let $w$ be a random word in $T$ of
  length $\ell$.  Then by (\ref{eq-abs-Z-mat}), the manifold $M_w$ has
  $|Z(M_w)| = \mu^{1-g} | \rho(w)_{(1,1)}|$.  We need to show that the
  distribution of $|Z(M_w)|$ limits on $\mu^{1-g} |X_{d_g}|$ as the length $\ell \to
  \infty$.

  Let $G$ be the closure of image of the subgroup $\langle T \rangle$ under
  $\rho$.  By Proposition~\ref{Prop:DenseImage} and
  Remark~\ref{rem-denseimage}, we know that $G$ contains $\SU{V}$.
  Now for the random word $w$, note that $\rho(w)$ is the product of the
  $\rho(t_i)$'s.  Thus by the Ito-Kawada Theorem \cite{ItoKawada1940},
  the matrix $\rho(w)$ becomes equidistributed on $G$ with respect to Haar
  measure as $\ell \to \infty$.  Thus if $X$ denotes the distribution of the
  $(1,1)$-entry of a random matrix in $G$, we have that $|Z(M)|$ is
  distributed by $\mu^{1-g} |X|$.  
  
  Thus it remains to show $X$ and $X_{d}$ have the same distribution
  where $d = d_g = \dim V$.  If $G = \U{d}$ this is obvious, so since
  $\SU{V} \leq G$ we must have $\dim G = \dim \SU{V} = \dim \U{V} -1$.
  Consider the homomorphism $\phi: S^1 \times G \to \U{d}$ given by $\phi(\zeta, M)
  = \zeta M$.  Since $S^1 \times G \to \U{d}$ is a covering map, the Haar
  measure (thought of as a differential form) on $\U{d}$ pulls back to
  (a multiple of) the Haar measure on $S^1 \times G$.  Now the pullback of
  the random variable $|X_{d}|$ on $\U{d}$ to $S^1 \times G$ is equal to
  $|X|$, and as the Haar measure pulls back we have that $|X_d|$ and
  $|X|$ have the same distribution.  To see that the distributions
  themselves are the same, for any $\zeta \in S^1$  consider left
  multiplication on $G$ by the matrix $A_\zeta \in G$ which is diagonal
  with entries $(\zeta, \zeta^{-1}, 1, 1, \ldots, 1)$.  Notice that the absolute
  value of the (1,1)-entry remains unchanged, but the phase shifts by
  $\zeta$.  Since Haar measure is left-invariant, this shows $X$ to be
  rotationally symmetric, and hence the same as $X_d$.
\end{proof}

We now consider the limit as the genus $g$ gets large.

\begin{proof}[Proof of Theorem~\ref{Thm:DistributionAbs}]
  We need to determine the limit of the random variables $\mu^{1-g}
  |X_{d_g}|$ as $g \to \infty$.  It is well known that $\sqrt{d_g} X_{d_g}$
  converges in distribution to the standard complex Gaussian (see
  e.g.~\cite{Novak2007}), and by Lemma~\ref{Lem:MuVersusD} we know
  $\sqrt{d_g} \approx \mu^{1-g}$.  Combining these can give that $\mu^{1-g}
  |X_{d_g}|$ converges to the standard Rayleigh distribution $R$, but
  instead we show this directly starting from
  Theorem~\ref{Lem:ZMdistributedasXd}.  If we fix $x \geq 0$, then since
  $\mu < 1$ we have
  \begin{eqnarray*}
  \lim_{g \to \infty}   P\left\{ |Z(M)| \geq x \right\} &=& \lim_{g \to \infty} \big(1 - \mu^{2 g - 2}x^2\big)^{d_g - 1}  \\
  &=&\lim_{g \to \infty} \left(1 - \frac{x^2}{\mu^{2 - 2g}} \right)^{\mu^{2 - 2g}} \cdot \lim_{g \to \infty} \left(1 - \frac{x^2}{\mu^{2 - 2g}} \right)^{\frac{d_g - 1}{\mu^{2 - 2g}}} \\
  &=& e^{-x^2} \cdot 1^1 = e^{-x^2} = P\left\{ R \geq x \right\}
    \end{eqnarray*}
    as needed.
\end{proof}

As $|Z(M)|$ is distributed by the standard Rayleigh distribution, it is very natural to postulate
\begin{conjecture} 
  \label{Conj} 
  Let $r \geq 5$ be prime.  For a random 3-manifold $M$, the $\SO{3}$
  invariant $Z(M)$ is distributed by a standard complex Gaussian.
\end{conjecture}
Looking at the proof of Theorem~\ref{Thm:DistributionAbs}, in fact we
showed that $\mu^{1-g} \rho(w)_{(1,1)}$ is essentially distributed by a
standard complex Gaussian when $g$ is large.  From (\ref{eq-Z-mat}),
we see this quantity differs from $Z(M)$ by $\kappa^{-3 (\sigma_b(w) +
  e(w))}$, where $\kappa$ is a certain $4r^{\mathrm{th}}$ root of unity.
Since $r$ is fixed, this tells us the image of $Z(M)$ in $\C /(z \mapsto
\kappa^3 z)$ is distributed by the push-forward of the standard complex
Gaussian.  The exponent sum $e(w)$ is a homomorphism $e \maps
\Free(\Dns) \to \Z$ and is thus easy to deal with, but $\sigma_b(w)$ is not
a homomorphism.  Indeed, the signature $\sigma_b$ is not Markovian in the
sense that $\sigma_b(w \cdot t)$ does not depend just on $\sigma_b(w)$ and $t$.
While it seems almost certain that $\sigma_b(w) + e(w)$ must be
uncorrelated with $\rho(w)_{(1,1)}$, and hence $Z(M)$ is the standard
complex Gaussian, we were unable to show this.  If
Conjecture~\ref{Conj} is true, it would give an alternate proof that
the values of $Z(M)$ are dense in $\C$ \cite{Wong:density}.

\begin{remark} A natural approach to Conjecture~\ref{Conj} is to
    work with one of the extended mapping class groups discussed in
    \cite{MasbaumRoberts1995}.  For instance, we used one
    such group $\MCGext$ above in the proof of
    Proposition~\ref{Prop:SU} since $\MCGext$ has an honest
    representation on $V$, not just a projective one.  Certainly, one
    can generate a random Heegaard splitting via a random walk in
    $\MCGext$ instead of $\MCG$, but one should think of the output as
    a \3-manifold $M$ \emph{and} a $p_1$-structure $\xi$.  That is,
    such a random model really produces \emph{$p_1$-manifolds} in the
    language of \cite[Section 4]{MasbaumRoberts1995}.  The $\SO{3}$ TQFT makes
    sense for such $p_1$-manifolds, and the plain $Z(M)$ is simply
    that more general invariant of $M$ and its canonical
    $p_1$-structure.  Now the above arguments show that $Z(M, \xi)$ is
    distribution by the standard complex Gaussian.  However, when one
    attempts to use this to understand the original question, the key
    issue is how a randomly chosen $\xi$ differs from the canonical
    one, which hinges on the signature issue discussed in the
    preceding paragraph.  Thus, while this line of thinking is a nice
    way to frame the underlying issue, it does not seem to give
    purchase on Conjecture~\ref{Conj}.
\end{remark}

Finally, we compute the probability that $Z(M)$ gives
a sharp genus bound for each genus $g$ and level $r$.

\begin{proof}[Proof of Theorem~\ref{thm-genus-bound}]
  Fix a genus $g$.  By \cite{Maher2008}, a random Heegaard splitting
  of genus $g$ actually has Heegaard genus $g$ with probability 1.
  Thus by Theorem~\ref{Thm:HGbound}, the quantum invariant will give a
  sharp lower bound on Heegaard genus whenever $|Z(M)| > \mu^{2-g}$.
  By Theorem~\ref{Thm:DistributionAbs}, this happens with probability
  $\left(1 - \mu^2 \right)^{d_g - 1}$, as claimed. 
\end{proof}

\section{Random surface bundles}
\label{sec-bundle}

A simple class of 3-manifolds are those which are surface bundles
fibering over the circle.  Such manifolds are built from some $f \in
\Gamma$ by considering the mapping torus
\[
MT_f = \rightquom{\Sigma \times [0,1]}{(p, 1) \sim (f(p), 0)}{3pt}{\Big}
\]
For reasons of homology, a random \3-manifold in the sense of
Section~\ref{sec-random} fibers over the circle with probability 0
\cite[Cor.~8.5]{DunfieldThurston:random}, and even among those
manifolds with $b_1 > 0$, surface bundles appear to be rare
\cite{DunfieldDThurston2006}.  Despite this, we show that for a
natural model of random bundles the distribution of $| Z(M) |$ is the
same as that of random manifolds more generally.  This can be taken as
evidence for the naturality and robustness of this notion of random
\3-manifolds.

Here, a random bundle with fiber a surface $\Sigma$ of genus $g$ is
defined just as one expects from Section~\ref{sec-random}: after
fixing generators $T$ of the mapping class group $\Gamma$ of $\Sigma$, one
considers $MT_f$ for $f \in \Gamma$ a random word in $T$ of
length $\ell$, as $\ell$ tends to infinity.   If we then send $g \to \infty$, we find:
\begin{theorem}\label{thm-bundle} 
  Consider the $\SO{3}$ invariant for a prime level $r \geq 5$.  Then
  $|Z(M)|$ is distributed by the standard Rayleigh for $M$ a random
  surface bundle.  
\end{theorem}
\begin{proof}
  Fix a genus $g > 2$ for the fiber $\Sigma$.  In the notation of
  Section~\ref{Section:CentralExt}, if $w \in \Free(\Dns)$ then we have
  (see e.g.~\cite[\S IV.7.2]{Turaev1994}):
  \[
  | Z(MT_w) | = \big| \tr \rho(w) \big| 
  \]
  Thus as in the proof of Theorem~\ref{Lem:ZMdistributedasXd}, the
  Ito-Kawada Theorem tells us that $| Z(MT_w) | $ is distributed the
  same as $|\tr A|$ for $A \in \SU{V}$ chosen with respect to Haar
  measure.  Again as in Theorem~\ref{Lem:ZMdistributedasXd}, this is
  the same as the distribution of $|\tr A|$ for $A \in \U{V}$.  Now as
  $d \to \infty$, the distribution of $\tr A$ for $A \in \U{d}$ converges to
  the standard complex Gaussian (\cite{DiaconisShahshahani1994}, see
  also \cite{PetzReffy2004, Novak2007}).  Thus since $\dim V \to \infty$ as
  $g \to \infty$ by Lemma~\ref{Lem:MuVersusD}, we have that $| Z(M) | $
  is distributed by the standard Rayleigh for a random surface bundle.
\end{proof}

A genus-$g$ surface bundle $M$ has a natural Heegaard splitting of
genus $2 g + 1$ obtained by tubing together two of the fibers.
Combining \cite{Maher2008} with Theorem~4.2 of \cite{Souto2007} shows
that this is the minimal genus splitting with probability 1, and
explicit examples where this is the case are easy to construct by
taking the monodromy from the Torelli group.  
However, we claim that
the $\SO{3}$ quantum invariants can never certify that $M$ has
Heegaard genus $2g +1$.  As in the proof of Theorem~\ref{thm-bundle},
we have $|Z(M)|= |\tr \rho(w)|$, and since $\rho(w)$ is unitary, it
follows that $| \tr \rho(w)| \leq \tr(\mathrm{Id})=d_g$.  From
(\ref{eq-dg-vs-mu}), we see that $d_g \leq \mu^{2-2g} < \mu^{1 - 2g}$ for
any $g$ when $r \geq 5$.  Thus $|Z(M)| < \mu^{1-2g},$ and so $Z(M)$ fails to give a
sharp bound on Heegaard genus.  To summarize, this shows:
\begin{proposition}
There exist closed hyperbolic \3-manifolds of arbitrary large
Heegaard genus such that no $\SO 3$ quantum invariant gives
a sharp genus bound.
\end{proposition}

{\RaggedRight \bibliographystyle{math} \bibliography{quantum} }

{\scriptsize \versioninfo}
\end{document}